\documentclass[onecolumn]{IEEEtran}

\usepackage{amssymb}
\usepackage{graphicx}
\usepackage{epsfig}
\usepackage{delarray}
\usepackage{amstext}

\newcommand{\U}[1]{\,^{{#1}\hspace*{-1pt}}}            

\title{Probabilistic DHP Adaptive Critic for Nonlinear Stochastic Control Systems}

\author{Randa~Herzallah
\thanks{R. Herzallah is with Al-Balqa' Applied University, Jordan.}}

\begin{document}

\maketitle

\begin{abstract}

Following the recently developed algorithms for fully probabilistic control design for general dynamic stochastic systems~\cite{herzallah11,Karny96}, this paper presents the solution to the probabilistic dual heuristic programming (DHP) adaptive critic method~\cite{herzallah11} and randomized control algorithm for stochastic nonlinear dynamical systems. The purpose of the randomized control input design is to make the joint probability density function of the closed loop system as close as possible to a predetermined ideal joint probability density function. This paper completes the previous work~\cite{herzallah11,Karny96} by formulating and solving the fully probabilistic control design problem on the more general case of nonlinear stochastic discrete time systems. A simulated example is used to demonstrate the use of the algorithm and encouraging results have been obtained.

\end{abstract}


\section{Introduction}\label{Introd}
The mean variance~\cite{Astrom70} and utility function in linear quadratic optimal control~\cite{Anderson71,Werbos92}, have been firstly introduced to characterize the performance of the closed loop stochastic control systems. In more recent work, stochastic adaptive control~\cite{Meyn87}, stochastic linear quadratic martingale~\cite{Solo90,Everdij96}, sliding mode control for stochastic systems~\cite{Yugang08}, and predictive stochastic control~\cite{Filatov95,Blackmore10} have been proposed. In most of previous works the system under control has been assumed to be linear or of Gaussian probability density function (pdf). However, it has been shown~\cite{Herzallah12,Wang09,Fabri01} that in systems where the stochastic signals is non Gaussian or the system dynamics have strong nonlinearity, existing control methods do not generally yield optimization of the control system.

Consequently, in the past few years four groups of control algorithms for general stochastic systems with inherent models' and parameters' uncertainty have been developed: (1) the control of the shape of the output pdf~\cite{Wang99,Wang01}, (2) the minimum entropy control~\cite{Wang02}, (3) the Bayesian techniques for modelling and control~\cite{Herzallah08,Trees07} and (4) the control of the closed loop pdf~\cite{Karny96,KarGuy:20060030}. The control objective in the first group is to find a control input which makes the shape of the measured output pdf follows a given desired distribution. The second group is a generalization of the minimum variance control for linear Gaussian systems. The entropy in this group of control algorithms is used to characterize the performance of the closed loop systems and the controller is designed such that the shape of the pdf of the closed loop system is made as narrow as possible. In the third group, a general class of stochastic estimation and control problems is formulated from the Bayesian decision-theoretic viewpoint which is shown to be a general framework to solve stochastic estimation and control problems. Motivated by the probabilistic description of the closed loop control system, the Kullback--Leibler divergence distance has been proposed in the fourth group of control as a performance measure rather than the mean variance. This method of control is known as fully probabilistic Design (FPD). For stochastic systems with measurable states $x_t$, the objective of the FPD is to determine the pdf of a randomized optimal control law, $u_t$ described by,
\begin{equation}
\label{eq:probcriticNonL1}
c(u_t \mid x_{t-1}),
\end{equation}
that minimizes the discrepancy between the actual joint pdf of the closed loop system, $f$ and an ideal joint pdf, $\U{I}f$ measured by the Kullback Leibler divergence distance,
\begin{equation}
\label{eq:probcriticNonL2}
\mathcal{D}\left(f||\U{I}f\right)\equiv \int
f(x_t,u_t)\ln\left(\frac{f(x_t,u_t)}{\U{I}f(x_t,u_t)}\right)\,\mathrm{d} x_t \mathrm{d} u_t.
\end{equation}
The FPD problem is further simplified by the assumption that all pdfs needed in the design paradigm are existent and known. The main advantage of the FPD is that it provides an explicit form of the randomized optimal controller. However, since the evaluation of the randomized optimal controller involves multivariate integration steps which need to be computed by backward recursion the problem renders to be nontrivial and computationally very intensive. To overcome the difficulties arising in the FPD, a probabilistic DHP adaptive critic method is proposed in~\cite{herzallah11}. The DHP adaptive critic method uses a critic network to circumvent the need for explicitly evaluating the optimal value function, therefore, dramatically reducing computational requirements. Further more, unlike the FPD method all pdfs in the probabilistic DHP adaptive critic approach are assumed to be unknown, therefore, estimated using recent development in neural networks. However, up to now, the previous methods of FPD~\cite{Karny96} and probabilistic DHP critic~\cite{herzallah11} are demonstrated only on linear stochastic Gaussian systems where the means of the associated density functions are restricted to be linear functions, hence the solution to the control problem is constrained to the linear Gaussian control theory. This is restrictive in many real world applications that are characterized by strong nonlinearity and uncertainty. In practical processes where the forward and inverse dynamics of the system have strong nonlinearities, the means of the associated density functions are in general going to be nonlinear functions. This motivates the work in the current paper, which is concerned with the development of the solution to the FPD problem~\cite{Karny96} for nonlinear stochastic systems where the means of the various density functions are allowed to be general nonlinear functions.

For this purpose we adopt the probabilistic DHP adaptive critic method proposed in~\cite{herzallah11}. Although systems under consideration are nonlinear stochastic systems, all pdfs are assumed to be Gaussians. This assumption can be shown to represent no real restriction provided that the conditional expectations of these pdfs are estimated using nonlinear models. However, it is worth mentioning that although the pdfs of the nonlinear stochastic systems are approximated by Gaussian density functions, it is in general difficult to integrate the nonlinear conditional expectations of these density functions to the probabilistic DHP adaptive critic method or even the FPD method. As such, radial basis function (RBF) neural network with Gaussian basis functions~\cite{Bishop95} is proposed in this paper to approximate the conditional expectations of the nonlinear stochastic models. An important property of such a neural network is that it forms a unifying link with density estimation. As will be clear from subsequent developments, the use of RBF networks with Gaussian basis functions to approximate the unknown nonlinear models facilitates the evaluation of the Gaussian integrations and integrates naturally to the framework of probabilistic DHP adaptive critic paradigm.

The main achievement of this paper is the solution of the nonlinear fully probabilistic optimal control problem for stochastic nonlinear systems. Although, the previous methods~\cite{herzallah11,Karny96} have discussed the general framework of nonlinear control systems, the problem rendered to be very hard to implement even under the Gaussian and linear models assumptions. Hence, the control solution to these methods are derived and demonstrated on linear systems only. By using the probabilistic DHP critic method proposed in~\cite{herzallah11} and RBF networks to estimate all required density functions, we develop and demonstrate the solution to the FPD control problem on stochastic nonlinear systems. The derived solution provides an efficient approach to the solution of the fully probabilistic control design for nonlinear stochastic systems.

To achieve the objective of this paper, it will be organized as follows. Section~\ref{ProbForm} formulates the problem and discusses the problem of estimating the pdf of the system dynamics. Section~\ref{SolProb} presents the probabilistic DHP adaptive critic solution to nonlinear stochastic systems. The control algorithm of nonlinear control problems based on the probabilistic DHP adaptive critic method is discussed in Section~\ref{Algo}. Section~\ref{SimExample} contains a simulation example to show the effectiveness of the proposed probabilistic controller. The conclusion is provided in Section~\ref{Conc}.

\section{Problem Formulation and Preliminaries}\label{ProbForm}
\subsection{Model Description and Control objective}
The system considered in this paper is a nonlinear stochastic dynamical control system described by the following general stochastic equation,
\begin{equation}
\label{eq:probcriticNonL7}
x_{t} = \tilde{h}(x_{t-1}) + \tilde{g}(x_{t-1})u_t + \tilde{\epsilon}_t,
\end{equation}
where $x_t \in \mathcal{R}^n$ is the measured state vector, $u_t \in \mathcal{R}^r$ is the control input vector, $\tilde{h}(x_{t-1}): \mathcal{R}^n \longmapsto \mathcal{R}^n$ and $\tilde{g}(x_{t-1}): \mathcal{R}^n \longmapsto \mathcal{R}^r$ are unknown nonlinear functions of the state, and $\tilde{\epsilon}_t \in \mathcal{R}^n$ is an additive noise vector. The control problem confronted here is to design a control strategy for the system in~(\ref{eq:probcriticNonL7}) to control the state of the system to a predefined desired state value.  However, because of the noise input $\tilde{\epsilon}_t$ the previous state and present and future controls do not completely specify the present state, but instead determine only the probability distribution of these states, $s(x_t \mid u_t, x_{t-1})$. It is assumed that the noise input $\tilde{\epsilon}_t$ is unknown and hence the probability distribution of the states is unknown.

Since only probability distribution of the states can be determined, the above objective of this control problem should be re-defined in terms of the probabilistic control theory. Therefore, to achieve this control objective we consider designing a probabilistic controller $c(u_t \mid x_{t-1})$ that shapes the joint pdf of the closed loop system, $f(x_t,u_t)$ and makes it as close as possible to a predefined desired pdf, $\U{I}f(x_t,u_t)$. This design method was originally presented in~\cite{Karny96} where the probabilistic controller is obtained such that it minimizes the Kullback–Leibler divergence distance defined in~(\ref{eq:probcriticNonL2}).  The minimum cost function resulting from minimization of~(\ref{eq:probcriticNonL2}) with respect to admissible control sequence $u_t$, $t \in \{1, \dots, H\}$, with $H$ being the control horizon, is then shown to be given by the following recurrence equation~\cite{herzallah11},
\begin{eqnarray}
\label{eq:probcriticNonL3} &&-\ln(\gamma(x_{t-1}))=
\min_{c(u_{t}|x_{t-1})} \int
s(x_{t}|u_{t},x_{t-1})c(u_{t}|x_{t-1}) \times\bigg[\underbrace{\ln\left(
\frac{s(x_{t}|u_{t},x_{t-1})c(u_{t}|x_{t-1})}
     {\U{I}s(x_{t}|u_{t},x_{t-1})\U{I}c(u_{t}|x_{t-1})}
\right)}_{\equiv \mbox{partial cost $\Longrightarrow U(x_{t},u_{t})$}} \nonumber \\
&&-\underbrace{\ln(\gamma(x_{t}))}_{\mbox{optimal
cost-to-go}}\bigg] \,\mathrm{d}(x_{t},u_{t}),
\end{eqnarray}
where $-\ln(\gamma(x_{t-1}))$ is the expected minimum cost--to--go function and
\begin{equation}
\label{eq:probcriticNonL4}
f(x_t,u_t) = s(x_{t}|u_{t},x_{t-1})c(u_{t}|x_{t-1}),
\end{equation}
is the decomposition of the actual joint pdf by the chain rule~\cite{Pet:81}, which represents the most complete probabilistic description of the closed loop system. Here the pdf $s(x_{t}|u_{t},x_{t-1})$ describes the dynamics of the observed state vector $x_t$. Similarly
\begin{equation}
\label{eq:probcriticNonL5}
\U{I}f(x_t,u_t) = \U{I}s(x_{t}|u_{t},x_{t-1})\U{I}c(u_{t}|x_{t-1}),
\end{equation}
is the decomposition of the ideal joint pdf of the closed loop system and $\U{I}s(x_{t}|u_{t},x_{t-1})$ and $\U{I}c(u_{t}|x_{t-1})$ represent the pdfs of the desired dynamics of the observed state vector and ideal controller respectively. The solution of the FPD is given in the following proposition.
\\
\textbf{Proposition 1:}
The pdf of optimal controller minimizing the cost--to--go function~(\ref{eq:probcriticNonL3}) is given by
\begin{eqnarray}
\label{eq:probcriticNonL6}
&&c^*(u_{t}|x_{t-1}) = \frac{\U{I}c(u_{t}|x_{t-1}) \exp[-\beta(u_{t},x_{t-1})]}{\gamma(x_{t-1})}, \nonumber \\
&&\gamma(x_{t-1}) = \int \U{I}c(u_{t}|x_{t-1}) \exp[-\beta(u_{t},x_{t-1})] \mathrm{d} u_t, \nonumber \\
&&\beta(u_{t},x_{t-1}) = \int s(x_{t}|u_{t},x_{t-1}) \bigg[\ln \frac{s(x_{t}|u_{t},x_{t-1})}{\U{I}s(x_{t}|u_{t},x_{t-1})} - \ln (\gamma(x_{t})) \bigg] \mathrm{d} x_t.
\end{eqnarray}

\textbf{\textit{Proof}:} This proposition can be proven by adapting the proof of Proposition $2$ in~\cite{KarGuy:20060030}.

It should be noted that although a closed form can be found for the pdf of the optimal controller, the multivariate integrations in~(\ref{eq:probcriticNonL6}) are only tractable for the linear Gaussian case, where the mean of the Gaussian distribution is linear in the state and control values. Besides even for the linear Gaussian case, the solution to the optimal control history need to be computed from~(\ref{eq:probcriticNonL6}) by backward recursion. This backward dynamic programming approach is computationally very expensive and grows exponentially with the dimensionality of the state vector. To avoid these difficulties of the FPD, a probabilistic DHP adaptive critic method is proposed in~\cite{herzallah11} to approximate the optimal cost-to-go function and the probabilistic controller. Unknown pdfs were also estimated using recent development from neural network models. However, numerical experiments
and previous analytical studies have demonstrated the usefulness of this control approach to obtain the control efforts only on linear stochastic Gaussian systems. The solution to the probabilistic DHP adaptive critic methods for the nonlinear stochastic systems~(\ref{eq:probcriticNonL7}) on the other hand was not discussed. The objective of this paper is to discuss the various steps to obtain this solution and demonstrate the theoretical development on a nonlinear stochastic simulation example of the form given in~(\ref{eq:probcriticNonL7}).

We first start by discussing the estimation problem of unknown probabilistic models of the stochastic system defined in~(\ref{eq:probcriticNonL7}) and reviewing the probabilistic DHP adaptive critic methods that will be needed for further development in the article.

\subsection{pdf of the system dynamics}\label{ForProbModel}
To estimate the probabilistic model of the nonlinear stochastic system~(\ref{eq:probcriticNonL7}) we adopt the method proposed in~\cite{herzallah11}, where neural network models are used to provide a prediction for the conditional expectation of the system state values and calculating the global average variance of its residual error. For such a system, there exists a neural network model~\cite{herzallah11} such that the inequality,
\begin{equation}
\label{eq:probcriticNonL8}
\mid x_{t} - N_{f}(u_{t}, x_{t-1}) \mid \le \delta,
\end{equation}
holds, where $\delta >0$ is a known small number and $N_{f}(u_{t}, x_{t-1}) = \hat{x_t}$ is a neural network approximation of the state $x_t$. Assuming a RBF neural network model of the form $N_{f}(u_{t}, x_{t-1}) = h(x_{t-1})+g(x_{t-1})u_t$ in which $h(x_{t-1})$ and $g(x_{t-1})$ are estimates of the nonlinear functions $\tilde{h}(x_{t-1})$ and $\tilde{g}(x_{t-1})$ respectively, the stochastic system~(\ref{eq:probcriticNonL7}) can be re-expressed as
\begin{equation}
x_{t} = h(x_{t-1})+g(x_{t-1})u_t + e(x_{t-1},u_t).
\end{equation}
Here, $e(x_{t-1},u_t)$ represents the approximation error satisfying $\mid e(x_{t-1},u_t) \mid \le \delta$. This means that the resulting conditional distribution of the system dynamics $s(x_t \mid u_t, x_{t-1})$ is Gaussian distribution function with conditional expectation of the distribution being given by the neural network approximation and a global average covariance given by~\cite{herzallah11},
\begin{equation}
\label{eq:probcriticNonL9}
\Sigma =  E\left((x_{t}-\hat{x}_{t})(x_{t}-\hat{x}_{t})^{T}\right),
\end{equation}
with $E(.)$ denoting the expected value.

\subsection{Review of the Probabilistic DHP Adaptive Critic Method}\label{RevProbCritic}
The probabilistic DHP adaptive critic method uses, two neural networks: an adaptive critic network that approximates the derivative of the optimal cost-to-go function~(\ref{eq:probcriticNonL3}) with respect to the state, $\lambda^*[x_{t-1}] = \partial [-\ln(\gamma(x_{t-1}))]/\partial x_{t-1}$ and an action network that produces optimal randomized control inputs, $u^*_t$. In the probabilistic DHP critic method, the optimal control law is computed by deriving~(\ref{eq:probcriticNonL3}) with respect to the control input~\cite{herzallah11},
\begin{eqnarray}
\label{eq:probcriticNonL10}
&&\frac{\partial [-\ln(\gamma(x_{t-1}))]}{\partial u_{t}}\bigg|_{u_{t}=u_{t}^*} = \int
s(x_{t}|u_{t},x_{t-1})c(u_{t}|x_{t-1}) \times\bigg[ \frac{\partial
U(x_{t},u_{t})} {\partial x_{t}} \frac{\partial x_{t}}{\partial u_{t}} +
\frac{\partial U(x_{t-1},u_{t})}{\partial u_{t}} \nonumber \\
&&+ \lambda[x_{t}]
\frac{\partial x_{t}}{\partial u_{t}} \bigg]
\mathrm{d}(x_{t},u_{t})=0.
\end{eqnarray}
So the action network is optimized such that the error between optimal control input $u_{t}^*$, obtained from~(\ref{eq:probcriticNonL10}) and estimated control input $u_{t}$ from the neural network is minimized. Once this network is optimized information about the error between optimal control $u_{t}^*$ and estimated control $u_{t}$ will become available. This allows estimation of the conditional distribution of the randomized controller $c(u_t \mid x_{t-1})$ which is assumed to be Gaussian with mean computed from the output of the controller network and a global covariance matrix computed from the residual error between the output of the controller network and the optimal control signal, $E\left((u_{t}^*-u_{t})(u_{t}^*-u_{t})^{T}\right)$.

Given estimation of control law from the controller network and the derivative of the output of the critic network $\lambda[x_{t}]$, the critic network is then optimized by computing its desired value as follows~\cite{herzallah11},
\begin{eqnarray}
\label{eq:probcriticNonL11}
&&\hspace*{-6ex}\lambda^*[x_{t-1}] = \int
s(x_{t}|u_{t},x_{t-1})c(u_{t}|x_{t-1})\bigg[ \frac{\partial
U(x_{t},u_{t})} {\partial x_{t}} \frac{\partial x_{t}}{\partial x_{t-1}}
+  \frac{\partial U(x_{t},u_{t})}{\partial x_{t}} \frac{\partial
x_{t}}{\partial u_{t}}
\frac{\partial u_{t}}{\partial x_{t-1}} + \frac{\partial
U(x_{t},u_{t})}{\partial u_{t}}\frac{\partial u_{t}}{\partial x_{t-1}} \nonumber \\
 &+&\lambda[x_{t}] \frac{\partial x_{t}}{\partial x_{t-1}} +
\lambda[x_{t}] \frac{\partial x_{t}}{\partial u_{t}} \frac{\partial
u_{t}}{\partial x_{t-1}} \bigg ] \mathrm{d}(x_{t},u_{t}).
\end{eqnarray}
The training process for the adaptive critic network is a two stage process. The training of the action network, which outputs the optimal control policy $u[x_{t}]$ and the training of the critic network, which approximates the derivative of the cost function $\lambda [x_{t-1}]$. As a first step in the training process, the critic and the action networks need to be designed and the initial weights of these networks should be randomized. Since the derivative of the partial cost function can be calculated, this in combination with the critic outputs and the system model derivatives, allows the use of~(\ref{eq:probcriticNonL11}) to calculate the target value of the critic, $\lambda^*[x_{t-1}]$. The difference between $\lambda^*[x_{t-1}]$ and the output of the critic, $\lambda[x_{t-1}]$ is used to correct the critic network, until it converges. The output from the converged critic is used in~(\ref{eq:probcriticNonL10}) solving for the target $u^*_{t}$, which is then used to correct the action network. This alternating process of training the action and the critic networks is repeated until an acceptable performance is reached.

\section{Solution to the Problem}\label{SolProb}

\subsection{Basic Elements}\label{ModelRep}
In this section we derive the probabilistic DHP adaptive critic solution to the nonlinear stochastic system defined in~(\ref{eq:probcriticNonL7}). For presentations clarity and simplicity the solution to this problem will be developed for a regulation problem where the objective is to reach a zero state with a spread determined by a specified covariance matrix. Generalization to a state value that is different than zero is straight forward.

As discussed in Section~\ref{ProbForm}, the conditional distribution of the nonlinear system~(\ref{eq:probcriticNonL7}) is estimated as a Gaussian distribution described by,
\begin{eqnarray}
\label{eq:probcriticNonL12}
x_{t} = h(x_{t-1})+g(x_{t-1})u_t + e(x_{t-1},u_t) \nonumber \\
s(x_t \mid u_t, x_{t-1}) \rightsquigarrow \mathcal{N}_{x_{t}}
(h(x_{t-1})+g(x_{t-1}) u_{t}, \Sigma).
\end{eqnarray}
For the considered regulation problem, the system is initially in state $x_{t-1}$ and the aim is to
return the system state to the origin. Hence, the distribution of the ideal state of the system
is taken to be,
\begin{equation}
\label{eq:probcriticNonL13}
\U{I}s(x_{t}|u_{t},x_{t-1}) = \mathcal{N}_{x_{t}}(0, \Sigma),
\end{equation}
where here the desired mean value of the state is taken to be zero and where $\Sigma$ specifies the covariance of the innovation of the state values.

The stochastic model of the randomized controller to be designed is estimated as discussed in Section~\ref{RevProbCritic} by the well known RBF neural network,
\begin{eqnarray}
\label{eq:probcriticNonL14}
u^k_{t} &=& \sum_{j=0}^M w_{kj} \psi_j(x_{t-1}) + \omega^k_{t} \nonumber \\
u_t &=& W \psi(x_{t-1}) + \omega_{t} \nonumber \\
c(u_{t}|x_{t-1}) & \rightsquigarrow & \mathcal{N}_{u_{t}}(W \psi(x_{t-1}), \Gamma),
\end{eqnarray}
where $W = [w_{kj}]$ is the matrix of the weight parameters, $M$ is the number of basis functions of the controller network, $\omega_{t}$ is the residual error of the control input vector, $\Gamma$ is the covariance of the residual error of control, and the RBF activation functions, $\psi_j(x_{t-1})$ are Gaussian basis functions~\cite{Bishop95},
\begin{equation}
\label{eq:probcriticNonL15}
\psi_j(x_{t-1}) = \exp \bigg ( -(x_{t-1} - \mu_j )^T \rho_j^{-1} (x_{t-1} - \mu_j ) \bigg).
\end{equation}
The distribution of the ideal controller is also assumed to be
\begin{equation}
\label{eq:probcriticNonL16}
\U{I}c(u_{t}|x_{t-1}) = \mathcal{N}_{u_{t}}(0, \Gamma).
\end{equation}
The desired value of the critic network given in~(\ref{eq:probcriticNonL11}) is also taken to be the target of an RBF neural network as follows~\cite{Bishop95},
\begin{eqnarray}
\label{eq:probcriticNonL17}
\lambda^m[x_{t-1}] &=& \sum_{l=0}^L \chi_{ml} \phi_l(x_{t-1}), \nonumber \\
\lambda[x_{t-1}] &=&  \chi \phi(x_{t-1}),
\end{eqnarray}
where $\chi$ is the matrix of weight parameters of the critic network, $L$ is the number of basis functions, and the basis functions, $\phi_l(x_{t-1})$ are taken to be Gaussian basis functions~\cite{Bishop95},
\begin{equation}
\label{eq:probcriticNonL18}
\phi_l(x_{t-1}) = \exp \bigg ( -(x_{t-1} - z_l)^T \gamma_l^{-1} (x_{t-1} - z_l) \bigg).
\end{equation}
Having defined those elements, the solution to the probabilistic DHP adaptive critic can now be obtained by calculating the desired value of the critic network and the optimal control inputs. We start in the next section by calculating the desired value of the critic. The optimal control input will be calculated in Section~\ref{ProbControl}.
\subsection{Desired Value of the Critic Network}\label{DesCritic}
For the conditional distribution of nonlinear system~(\ref{eq:probcriticNonL12}), conditional distribution of nonlinear controller~(\ref{eq:probcriticNonL14}), and nonlinear critic model~(\ref{eq:probcriticNonL17}), the desired target value of the critic network can be calculated by carrying out the calculations implied by~(\ref{eq:probcriticNonL11}). Starting by the first term on the right hand side of (\ref{eq:probcriticNonL11}) we get,
\begin{eqnarray}
\label{eq:probcriticNonL19}
&&\int s(x_{t}|u_{t},x_{t-1})c(u_{t}|x_{t-1}) \frac{\partial
U(x_{t},u_{t})} {\partial x_{t}} \frac{\partial x_{t}}{\partial x_{t-1}}
\mathrm{d}(x_{t},u_{t})= \int \exp[-(u_t-\hat{u}_t)^T \Gamma^{-1} (u_t-\hat{u}_t)] \times  \nonumber \\
&&\bigg\{ \int \exp[-(x_t-\hat{x}_t)^T \Sigma^{-1} (x_t-\hat{x}_t)] 2 \hat{x}_t^T \Sigma^{-1} [ h'(x_{t-1}) + g'(x_{t-1}) u_t] \mathrm{d} x_t \bigg \} \mathrm{d} u_t \nonumber \\ && = \int \exp[-(u_t-\hat{u}_t)^T \Gamma^{-1} (u_t-\hat{u}_t)] 2 \hat{x}_t^T \Sigma^{-1} [ h'(x_{t-1}) + g'(x_{t-1}) u_t] \mathrm{d} u_t \nonumber \\ && = 2 [h(x_{t-1}) + g(x_{t-1}) \hat{u}_t]^T \Sigma^{-1} [h'(x_{t-1}) + g'(x_{t-1}) \hat{u}_t],
\end{eqnarray}
where we have introduced the definitions $\hat{x}_t = h(x_{t-1}) + g(x_{t-1}) u_t$ and $\hat{u}_t = W \psi(x_{t-1})$ and where $h'(x_{t-1})= \frac{\partial h(x_{t-1})}{\partial x_{t-1}}$ and $g'(x_{t-1})= \frac{\partial g(x_{t-1})}{\partial x_{t-1}}$. The definition of the partial cost $U(x_{t},u_{t})$ is given in~(\ref{eq:probcriticNonL3}). For the considered regularization problem, it evaluates to $U(x_{t},u_{t}) = 2 \hat{x}_t^T \Sigma^{-1} x_t - \hat{x}_t^T \Sigma^{-1} \hat{x}_t + 2 \hat{u}_t^T \Gamma^{-1} u_t - \hat{u}_t^T \Gamma^{-1} \hat{u}_t$.

For the second term the partial derivatives of $U(x_{t},u_{t})$,
$x_{t}$, and $u_{t}$ with respect to $x_{t}$, $u_{t}$, and $x_{t-1}$
respectively need to be calculated,
\begin{eqnarray}
\label{eq:probcriticNonL20}
&& \int s(x_{t}|u_{t},x_{t-1})c(u_{t}|x_{t-1}) \frac{\partial
U(x_{t},u_{t})} {\partial x_{t}} \frac{\partial x_{t}}{\partial u_{t}}
\frac{\partial u_{t}}{\partial x_{t-1}} \mathrm{d}(x_{t},u_{t}) = \int \exp[-(u_t-\hat{u}_t)^T \Gamma^{-1} (u_t-\hat{u}_t)] \times  \nonumber \\ &&\bigg\{ \int \exp[-(x_t-\hat{x}_t)^T \Sigma^{-1} (x_t-\hat{x}_t)] 2 \hat{x}_t^T \Sigma^{-1} g(x_{t-1}) W \psi'(x_{t-1}) \mathrm{d} x_t \bigg \} \mathrm{d} u_t \nonumber \\ && = \int \exp[-(u_t-\hat{u}_t)^T \Gamma^{-1} (u_t-\hat{u}_t)] 2 \hat{x}_t^T \Sigma^{-1} g(x_{t-1}) W \psi'(x_{t-1}) \mathrm{d} u_t \nonumber \\ && = 2 [h(x_{t-1}) + g(x_{t-1}) \hat{u}_t]^T \Sigma^{-1} g(x_{t-1}) W \psi'(x_{t-1}),
\end{eqnarray}
where $\psi'(x_{t-1})= \frac{\partial \psi(x_{t-1})}{\partial x_{t-1}}$.

The third term requires calculation of the partial derivatives of $U(x_{t},u_{t})$ and $u_{t}$ with respect to $u_{t}$ and $x_{t-1}$ respectively,
\begin{eqnarray}
\label{eq:probcriticNonL21}
&&\int s(x_{t}|u_{t},x_{t-1})c(u_{t}|x_{t-1}) \frac{\partial
U(x_{t},u_{t})} {\partial u_{t}} \frac{\partial u_{t}}{\partial x_{t-1}}
\mathrm{d}(x_{t},u_{t}) = \int \exp[-(u_t-\hat{u}_t)^T \Gamma^{-1} (u_t-\hat{u}_t)] \times  \nonumber \\
&&\bigg\{ \int \exp[-(x_t-\hat{x}_t)^T \Sigma^{-1} (x_t-\hat{x}_t)] 2 \hat{u}_t^T \Gamma^{-1} W \psi'(x_{t-1}) \mathrm{d} x_t \bigg \} \mathrm{d} u_t \nonumber \\ && = \int \exp[-(u_t-\hat{u}_t)^T \Gamma^{-1} (u_t-\hat{u}_t)] 2 \hat{u}_t^T \Gamma^{-1} W \psi'(x_{t-1}) \mathrm{d} u_t \nonumber \\ && = 2 \hat{u}_t^T \Gamma^{-1} W \psi'(x_{t-1}).
\end{eqnarray}
The propagation of $\lambda[x_{t}]$ through the
stochastic model of (\ref{eq:probcriticNonL12}) back to $x_{t}$ yields the fourth term
\begin{eqnarray}
\label{eq:probcriticNonL22}
&&\int s(x_{t}|u_{t},x_{t-1})c(u_{t}|x_{t-1}) \lambda[x_{t}]
\frac{\partial x_{t}}{\partial x_{t-1}} \mathrm{d}(x_{t},u_{t})=\int \exp[-(u_t-\hat{u}_t)^T \Gamma^{-1} (u_t-\hat{u}_t)] \times  \nonumber \\ && \alpha(u_t, x_{t-1}) \mathrm{d} u_t,
\end{eqnarray}
where we used
\begin{eqnarray}
\label{eq:probcriticNonL23}
&&\alpha(u_t, x_{t-1}) = \int \exp[-(x_t-\hat{x}_t)^T \Sigma^{-1} (x_t-\hat{x}_t)] (\chi \phi(x_{t}))^T [h'(x_{t-1}) + g'(x_{t-1}) u_t] \mathrm{d} x_t,~~~~~~
\end{eqnarray}
and where we used~(\ref{eq:probcriticNonL17}) with $x_t$ as input to obtain $\lambda[x_{t}]=\chi \phi(x_{t})$. Using~(\ref{eq:probcriticNonL18}) but again with $x_t$ as an input in~(\ref{eq:probcriticNonL23}) yields,
\begin{eqnarray}
\label{eq:probcriticNonL24}
&&\alpha(u_t, x_{t-1}) = \int \left [ \begin{array}{c} \exp[-(x_t-\hat{x}_t)^T \Sigma^{-1} (x_t-\hat{x}_t)] \exp[-(x_t-z_1)^T \gamma_1^{-1} (x_t-z_1)] \\ \exp[-(x_t-\hat{x}_t)^T \Sigma^{-1} (x_t-\hat{x}_t)] \exp[-(x_t-z_2)^T \gamma_2^{-1} (x_t-z_2)] \\ \vdots \\ \exp[-(x_t-\hat{x}_t)^T \Sigma^{-1} (x_t-\hat{x}_t)] \exp[-(x_t-z_L)^T \gamma_L^{-1} (x_t-z_L)]  \end{array} \right ]^T \chi^T \nonumber \\&& \times [h'(x_{t-1}) + g'(x_{t-1}) u_t]  \mathrm{d} x_t  \\ \label{eq:probcriticNonL24'} && = \left [ \begin{array}{c} \exp[-(\hat{x}_t-z_1)^T (\Sigma+\gamma_1)^{-1} (\hat{x}_t-z_1)] \\ \exp[-(\hat{x}_t-z_2)^T (\Sigma+\gamma_2)^{-1} (\hat{x}_t-z_2)] \\ \vdots \\ \exp[-(\hat{x}_t-z_L)^T (\Sigma+\gamma_L)^{-1} (\hat{x}_t-z_L)]  \end{array} \right ]^T \chi^T [h'(x_{t-1}) + g'(x_{t-1}) u_t],
\end{eqnarray}
where the sum of two quadratics as a result of the multiplication of the two exponentials of the first term in brackets in~(\ref{eq:probcriticNonL24}) is rewritten in the following form $(x_t-\hat{x}_t)^T \Sigma^{-1} (x_t-\hat{x}_t) + (x_t-z_l)^T \gamma_l^{-1} (x_t-z_l) = (x_t-{\bf E}_l)^T (\Sigma^{-1}+\gamma^{-1}_l) (x_t-{\bf E}_l) + (\hat{x}_t-z_l)^T (\Sigma+\gamma_l)^{-1} (\hat{x}_t-z_l)$ and the integration with respect to $x_t$ is evaluated. Here ${\bf E}_l = (\Sigma^{-1}+\gamma^{-1}_l)^{-1}(\Sigma^{-1} \hat{x}_t + \gamma^{-1}_l z_l)$.  The substitution of~(\ref{eq:probcriticNonL24'}) in~(\ref{eq:probcriticNonL22}) yields the value of the fourth term
\begin{eqnarray}
\label{eq:probcriticNonL25}
&&\int s(x_{t}|u_{t},x_{t-1})c(u_{t}|x_{t-1}) \lambda[x_{t}]
\frac{\partial x_{t}}{\partial x_{t-1}} \mathrm{d}(x_{t},u_{t})=\int \exp[-(u_t-\hat{u}_t)^T \Gamma^{-1} (u_t-\hat{u}_t)] \times  \nonumber \\ && \left [ \begin{array}{c} \exp[-(\hat{x}_t-z_1)^T (\Sigma+\gamma_1)^{-1} (\hat{x}_t-z_1)] \\ \exp[-(\hat{x}_t-z_2)^T (\Sigma+\gamma_2)^{-1} (\hat{x}_t-z_2)] \\ \vdots \\ \exp[-(\hat{x}_t-z_L)^T (\Sigma+\gamma_L)^{-1} (\hat{x}_t-z_L)]  \end{array} \right ]^T \chi^T [h'(x_{t-1}) + g'(x_{t-1}) u_t] \mathrm{d} u_t \nonumber \\ && = \int \left [ \begin{array}{c} \exp[-(\hat{x}_t-z_1)^T (\Sigma+\gamma_1)^{-1} (\hat{x}_t-z_1) - (u_t-\hat{u}_t)^T \Gamma^{-1} (u_t-\hat{u}_t)] \\ \exp[-(\hat{x}_t-z_2)^T (\Sigma+\gamma_2)^{-1} (\hat{x}_t-z_2)- (u_t-\hat{u}_t)^T \Gamma^{-1} (u_t-\hat{u}_t)] \\ \vdots \\ \exp[-(\hat{x}_t-z_L)^T (\Sigma+\gamma_L)^{-1} (\hat{x}_t-z_L)- (u_t-\hat{u}_t)^T \Gamma^{-1} (u_t-\hat{u}_t)]  \end{array} \right ]^T \chi^T \nonumber \\&& \times [h'(x_{t-1}) + g'(x_{t-1}) u_t] \mathrm{d} u_t.
\end{eqnarray}
To simplify the integration in~(\ref{eq:probcriticNonL25}), the square of the argument of the exponential is completed,
\begin{eqnarray}
\label{eq:probcriticNonL26}
&&(\hat{x}_t-z_l)^T (\Sigma+\gamma_l)^{-1} (\hat{x}_t-z_l)- (u_t-\hat{u}_t)^T \Gamma^{-1} (u_t-\hat{u}_t) = \nonumber \\
&&(u_t - H_l)^T \Omega_l (u_t - H_l) + (z_l - h(x_{t-1}))^T (\Sigma+\gamma_l)^{-1}(z_l - h(x_{t-1})) \nonumber \\ &&+ \hat{u}_t^T \Gamma^{-1} \hat{u}_t - H_l^T \Omega_l H_l,~~~~~~~~~
\end{eqnarray}
where $\Omega_l = g^T(x_{t-1}) (\Sigma + \gamma_l)^{-1} g(x_{t-1}) + \Gamma^{-1}$, and $H_l = \Omega_l^{-1}[g^T(x_{t-1}) (\Sigma + \gamma_l)^{-1}(z_l - h(x_{t-1})) + \Gamma^{-1} \hat{u}_t]$. Denoting the exponential of the last constant three terms (the terms that are independent of $u_t$) on the right hand side of~(\ref{eq:probcriticNonL26}) by $\digamma_l = \exp[-(z_l - h(x_{t-1}))^T (\Sigma+\gamma_l)^{-1}(z_l - h(x_{t-1})) - \hat{u}_t^T \Gamma^{-1} \hat{u}_t + H_l^T \Omega_l H_l]$ and substituting~(\ref{eq:probcriticNonL26}) in~(\ref{eq:probcriticNonL25}) yields,
\begin{eqnarray}
\label{eq:probcriticNonL27}
&&\int s(x_{t}|u_{t},x_{t-1})c(u_{t}|x_{t-1}) \lambda[x_{t}]
\frac{\partial x_{t}}{\partial x_{t-1}} \mathrm{d}(x_{t},u_{t})= \nonumber \\ && \int  \left [ \begin{array}{c} \exp[-(u_t-H_1)^T \Omega_1 (u_t-H_1)] \\ \exp[-(u_t-H_2)^T \Omega_2 (u_t-H_2)] \\ \vdots \\ \exp[-(u_t-H_L)^T \Omega_L (u_t-H_L)]  \end{array} \right ]^T. \digamma^T \chi^T  [h'(x_{t-1}) + g'(x_{t-1}) u_t] \mathrm{d} u_t \nonumber \\&& = [\chi \digamma]^T h'(x_{t-1}) +  [\chi \digamma]^T g'(x_{t-1}) H.
\end{eqnarray}
Finally the fifth term can be calculated by propagating
$\lambda[x_{t}]$ through the stochastic model of
(\ref{eq:probcriticNonL12}) back to $u_{t}$, and then through the action
network, which yields,
\begin{eqnarray}
\label{eq:probcriticNonL28}
\hspace*{-2ex}
&&\int s(x_{t}|u_{t},x_{t-1})c(u_{t}|x_{t-1}) \lambda[x_{t}]
\frac{\partial x_{t}}{\partial u_{t}} \frac{\partial u_{t}}{\partial x_{t-1}}
\mathrm{d}(x_{t},u_{t})= \int \exp[-(u_t-\hat{u}_t)^T \Gamma^{-1} (u_t-\hat{u}_t)] \times  \nonumber \\
&&\bigg\{ \int \exp[-(x_t-\hat{x}_t)^T \Sigma^{-1} (x_t-\hat{x}_t)] [\chi \phi(x_{t})]^T g(x_{t-1}) W \psi'(x_{t-1}) \mathrm{d} x_t \bigg \} \mathrm{d} u_t \nonumber \\ && = [\chi \digamma]^T g(x_{t-1}) W \psi'(x_{t-1}).
\end{eqnarray}
Adding~(\ref{eq:probcriticNonL19}),~(\ref{eq:probcriticNonL20}),~(\ref{eq:probcriticNonL21}),~(\ref{eq:probcriticNonL27}),~(\ref{eq:probcriticNonL28}) together, yields the target vector of the critic network,
\begin{eqnarray}
\label{eq:probcriticNonL29}
\lambda^*[x_{t-1}] &=& 2 [h(x_{t-1}) + g(x_{t-1}) \hat{u}_t]^T \Sigma^{-1} [h'(x_{t-1}) + g'(x_{t-1}) \hat{u}_t] \nonumber \\ &+& 2 [h(x_{t-1}) + g(x_{t-1}) \hat{u}_t]^T \Sigma^{-1} g(x_{t-1}) W \psi'(x_{t-1}) + 2 \hat{u}_t^T \Gamma^{-1} W \psi'(x_{t-1}) \nonumber \\ &+& [\chi \digamma]^T h'(x_{t-1})  +  [\chi \digamma]^T g'(x_{t-1}) H + [\chi \digamma]^T g(x_{t-1}) W \psi'(x_{t-1}).
\end{eqnarray}
Equation~(\ref{eq:probcriticNonL29}) can then be used to correct the critic network and update its parameters.

\textit{Remark:} All of the above integrals (implied by~(\ref{eq:probcriticNonL19})--~(\ref{eq:probcriticNonL28})) are Gaussian integrals which can be evaluated using theorems and corollaries provided in~\cite{Franklin83}, chapter $10$.

\subsection{Probabilistic Control}\label{ProbControl}
The randomized control input can be computed by solving the optimality equation~(\ref{eq:probcriticNonL10}). The first term on the right hand side of~(\ref{eq:probcriticNonL10}) is the expected value of the partial derivatives of $U(x_{t},u_{t})$ and $x_{t}$ with respect to $x_{t}$ and $u_{t}$ respectively,
\begin{eqnarray}
\label{eq:probcriticNonL30}
&&\int s(x_{t}|u_{t},x_{t-1})c(u_{t}|x_{t-1}) \frac{\partial
U(x_{t},u_{t})} {\partial x_{t}} \frac{\partial x_{t}}{\partial u_{t}}
\mathrm{d}(x_{t},u_{t}) = \int \exp[-(u_t-\hat{u}_t)^T \Gamma^{-1} (u_t-\hat{u}_t)] \times  \nonumber \\
&&\bigg\{ \int \exp[-(x_t-\hat{x}_t)^T \Sigma^{-1} (x_t-\hat{x}_t)] 2 \hat{x}_t^T \Sigma^{-1} g(x_{t-1}) \mathrm{d} x_t \bigg \} \mathrm{d} u_t \nonumber \\ && = \int \exp[-(u_t-\hat{u}_t)^T \Gamma^{-1} (u_t-\hat{u}_t)] 2 [ h(x_{t-1}) + g(x_{t-1}) u_t]^T \Sigma^{-1} g(x_{t-1})  \mathrm{d} u_t \nonumber \\ && = 2 [h(x_{t-1}) + g(x_{t-1}) \hat{u}_t]^T \Sigma^{-1} g(x_{t-1}).
\end{eqnarray}
The second term requires the evaluation of the expected value of the partial derivatives of $U(x_{t},u_{t})$ with respect to $u_{t}$,
\begin{eqnarray}
\label{eq:probcriticNonL31}
&&\int s(x_{t}|u_{t},x_{t-1})c(u_{t}|x_{t-1}) \frac{\partial
U(x_{t},u_{t})} {\partial u_{t}} \mathrm{d}(x_{t},u_{t})= \int \exp[-(u_t-\hat{u}_t)^T \Gamma^{-1} (u_t-\hat{u}_t)] \times  \nonumber \\ &&\bigg\{ \int \exp[-(x_t-\hat{x}_t)^T \Sigma^{-1} (x_t-\hat{x}_t)] 2 \hat{u}_t^T \Gamma^{-1} \mathrm{d} x_t \bigg \} \mathrm{d} u_t \nonumber \\ && = \int \exp[-(u_t-\hat{u}_t)^T \Gamma^{-1} (u_t-\hat{u}_t)] 2 \hat{u}_t^T \Gamma^{-1} \mathrm{d} u_t \nonumber \\ && = 2 \hat{u}_t^T \Gamma^{-1}.
\end{eqnarray}
The last term is the expected value of the propagation of $\lambda[x_{t}]$ through the stochastic model of~(\ref{eq:probcriticNonL12}) back to $u_{t}$,
\begin{eqnarray}
\label{eq:probcriticNonL32}
&&\int s(x_{t}|u_{t},x_{t-1})c(u_{t}|x_{t-1}) \lambda[x_{t}] \frac{\partial
x_{t}} {\partial u_{t}} \mathrm{d}(x_{t},u_{t})= \int \exp[-(u_t-\hat{u}_t)^T \Gamma^{-1} (u_t-\hat{u}_t)] \times  \nonumber \\ &&\bigg\{ \int \exp[-(x_t-\hat{x}_t)^T \Sigma^{-1} (x_t-\hat{x}_t)] (\chi \phi(x_t))^T g(x_{t-1}) \mathrm{d} x_t \bigg \} \mathrm{d} u_t \nonumber \\ && = \int \exp[-(u_t-\hat{u}_t)^T \Gamma^{-1} (u_t-\hat{u}_t)] \times  \nonumber \\ &&\bigg\{ \int \left [ \begin{array}{c} \exp[-(x_t-\hat{x}_t)^T \Sigma^{-1} (x_t-\hat{x}_t)] \exp[-(x_t-z_1)^T \gamma_1^{-1} (x_t-z_1)] \\ \exp[-(x_t-\hat{x}_t)^T \Sigma^{-1} (x_t-\hat{x}_t)] \exp[-(x_t-z_2)^T \gamma_2^{-1} (x_t-z_2)] \\ \vdots \\ \exp[-(x_t-\hat{x}_t)^T \Sigma^{-1} (x_t-\hat{x}_t)] \exp[-(x_t-z_L)^T \gamma_L^{-1} (x_t-z_L)]  \end{array} \right ]^T \chi^T g(x_{t-1}) \mathrm{d} x_t \bigg \} \mathrm{d} u_t \nonumber \\ && = \int \left [ \begin{array}{c} \exp[-(\hat{x}_t-z_1)^T (\Sigma+\gamma_1)^{-1} (\hat{x}_t-z_1) - (u_t-\hat{u}_t)^T \Gamma^{-1} (u_t-\hat{u}_t)] \\ \exp[-(\hat{x}_t-z_2)^T (\Sigma+\gamma_2)^{-1} (\hat{x}_t-z_2)- (u_t-\hat{u}_t)^T \Gamma^{-1} (u_t-\hat{u}_t)] \\ \vdots \\ \exp[-(\hat{x}_t-z_L)^T (\Sigma+\gamma_L)^{-1} (\hat{x}_t-z_L)- (u_t-\hat{u}_t)^T \Gamma^{-1} (u_t-\hat{u}_t)]  \end{array} \right ]^T \chi^T  \nonumber \\ && \times g(x_{t-1}) \mathrm{d} u_t \nonumber \\ && = (\chi \digamma)^T g(x_{t-1}).
\end{eqnarray}
Adding all terms together yields,
\begin{equation}
\label{eq:probcriticNonL33}
2 [h(x_{t-1}) + g(x_{t-1}) \hat{u}_t]^T \Sigma^{-1} g(x_{t-1}) + 2 \hat{u}_t^T \Gamma^{-1} + (\chi \digamma)^T g(x_{t-1}) = 0.
\end{equation}
The solution of this equation cannot be analytically obtained due to the nonlinear nature of $\digamma$. Hence the controller design is a nonlinear optimization problem which generally leads to a numerical solution. The controller network can then be optimized such that the error between optimal control input $u^*_t$ obtained from~(\ref{eq:probcriticNonL33}) and estimated control input $u_t$ from the neural network is minimized. The controller can then generates control signals $u_{t}$ stochastically from a gaussian distribution having a mean $\hat{u}_t$ and a global average covariance $\Gamma$ calculated as discussed in Section~\ref{RevProbCritic}.
\section{Nonlinear Randomized control Algorithm Based on Probabilistic DHP Critic Methods}\label{Algo}
The exact evaluation of the closed form of optimal controller in FPD methods~(\ref{eq:probcriticNonL6}) is nontrivial and computationally very intensive due to its involvement of multivariate integrations. These multivariate integrations are only tractable for the linear Gaussian case where the mean is linear in both state and control values. This motivated the probabilistic DHP adaptive critic approach discussed in Section~\ref{SolProb}. As can be seen from~(\ref{eq:probcriticNonL29}), the desired critic value can be calculated from neural network models of the forward dynamics of the system, nonlinear controller, and critic network. Similarly, the control law can be calculated from~(\ref{eq:probcriticNonL33}) once the critic network and system dynamic models become available. Although the probabilistic DHP adaptive critic approach involves multiple computation levels, its implementation can be made by means of modular approach constituting of functional modules and algorithmic modules~\cite{Silvia04,George02,Herzallah07c}. The key functional modules are the action and critic networks. Algorithmic modules on the other hand include the computation of the desired critic value~(\ref{eq:probcriticNonL29}), computation of the optimal control law~(\ref{eq:probcriticNonL33}), and the update of the networks parameters. Each module in the adaptive critic can be modified independently from other modules, thus facilitate fast and reliable implementation. The probabilistic DHP adaptive critic approach allows computation of probabilistic control laws from the FPD methods for linear Gaussian systems with less computational complexities, and more important allows the solution to nonlinear and non Gaussian systems.

\begin{figure}[htbp]
  \centerline{
    {\scalebox{0.5}[0.5]{\includegraphics{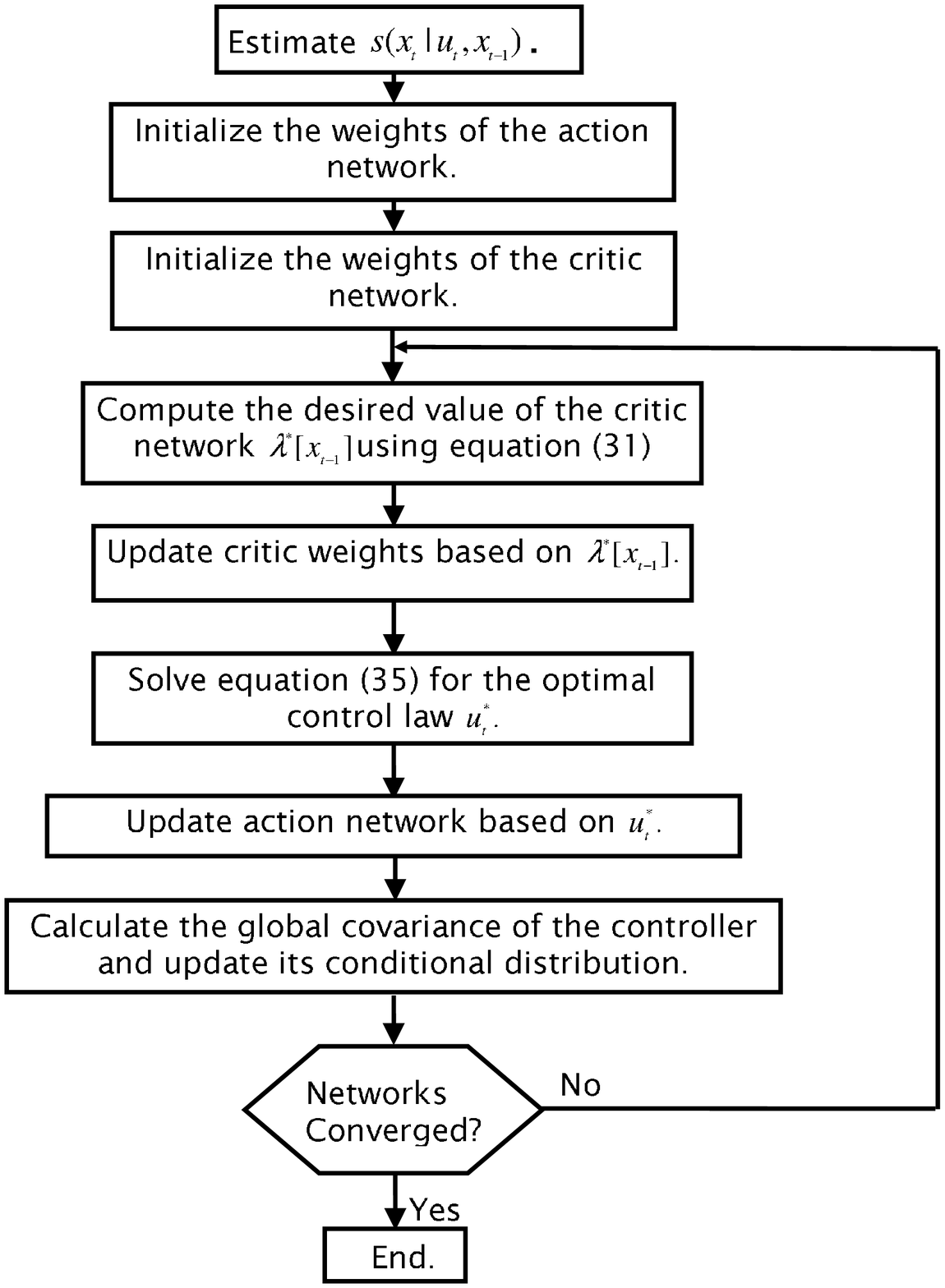}}
      }
    }
\caption[Implementation of the probabilistic DHP adaptive critic for nonlinear systems]
{Implementation of the probabilistic DHP adaptive critic for nonlinear systems.
{\label{fig:Flowchart1}
   }
  }
\end{figure}

The description below is appropriate for direct application to nonlinear control problems of the form stated in Section~\ref{ProbForm}.
\begin{enumerate}
\item Estimate the pdf of the stochastic model described by~(\ref{eq:probcriticNonL7}) as discussed in Section~\ref{ProbForm}.
\item Design and initialize the weights of the action network.
\item Design and initialize the weights of the critic network.
\item Calculate the desired value of the critic network using equation~(\ref{eq:probcriticNonL29}).
\item Use the difference between the desired value of the critic network as calculated from the previous step and the output of the critic to update the parameters of the critic network until it converges.
\item Use the output of the converged critic in~(\ref{eq:probcriticNonL33}) and solve for the optimal control law.
\item Use this value to update the parameters of the action network.
\item Calculate the global covariance matrix from the residual error between the output of the action network and optimal control law as calculated in step $6$ and update the conditional distribution of the optimal controller.
\item Repeat Steps $4-6$ until an acceptable performance is reached.
\end{enumerate}
The flow chart of implementing the probabilistic DHP adaptive critic method for nonlinear systems is given in Figure~\ref{fig:Flowchart1}. To summarize, the probabilistic DHP adaptive critic training algorithm cycles between a policy improvement routine and a control input determination operation, where the optimal control law and the derivative of the value function, $\lambda[x_{t-1}]$ are approximated by the action and critic networks respectively. The algorithm terminates when control law and the derivative of the value function have converged to the optimal or suboptimal control law and derivative value function respectively. The proof of convergence given in~\cite{Howard60,Silvia04} is directly applicable to the probabilistic adaptive critic design in this paper. A nonlinear control problem example to demonstrate the convergence of the proposed probabilistic critic network is given in Section~\ref{SimExample}. Further discussion on the convergence and speed of convergence of adaptive critic designs can be found in~\cite{George02}. Moreover, empirical evidence on the convergence of the adaptive critic design can be found in~\cite{Balakrishnan96,Balakrishnan02,Nilesh03,Chuan11}.

\section{Simulation Example}\label{SimExample}
This section demonstrates the probabilistic DHP adaptive critic method on nonlinear stochastic control system. The nonlinear dynamical system is described by the following equation
\begin{equation}
\label{eq:probcriticNonL34}
x_{t}= sin(x_{t-1}) + cos(3*x_{t-1}) + (2+ cos(x_{t-1})) u_t + \varepsilon_{t}.
\end{equation}
The unknown nonlinear dynamics are given by,
\begin{eqnarray}
\tilde{h}(x_{t-1}) &=& sin(x_{t-1}) + cos(3*x_{t-1}), \nonumber \\
\tilde{g}(x_{t-1}) &=& (2+ cos(x_{t-1})). \nonumber
\end{eqnarray}
The noise $\varepsilon_t$ is assumed to be sampled from a Gaussian distribution, $N(0,0.01)$. This system has been used in~\cite{Fabri01} to illustrate theoretical developments for suboptimal dual adaptive control.

The plant is initially, at time $t=0$, in state $x_0$, and the aim is to return the plant state to the origin, or a state close to the origin. Thus, a probabilistic DHP critic network that minimizes the cost function~(\ref{eq:probcriticNonL3}) is used to design the optimal probabilistic controller and derive optimal control inputs. As a first step in the solution the stochastic model of the forward dynamics of the plant described by~(\ref{eq:probcriticNonL34}) is estimated by a Gaussian density function as discussed in Section~\ref{ForProbModel}. Two RBF networks with $15$ and $6$ Gaussian basis functions are used to estimate the two nonlinear functions $\tilde{h}(x_{t-1})$ and $\tilde{g}(x_{t-1})$ respectively. Hence, the mean of the forward probabilistic model of Gaussian density function is given by $\hat{x}_t = h(x_{t-1}) + g(x_{t-1}) u_t$  and the its global variance, $\sigma^2 = 0.0098$ is computed from the residual error of the system dynamics. The weight parameters of the forward Gaussian probabilistic model are then kept fixed during the critic and action network training. To achieve the control objective of attaining a zero state, the distribution of the ideal state of the system dynamics is taken to be $\U{I}s(x_{t}|u_{t},x_{t-1}) = \mathcal{N}_{x_{t}}(0, 0.0098)$. Similarly, the ideal controller distribution is assumed to be $\U{I}c(u_{t}|x_{t-1}) = \mathcal{N}_{u_{t}}(0, 0.01)$.

The control is then initiated by another RBF network with six neurons in the hidden layer for random values of $x_t$, taken uniformly from the interval $[-4,4]$. Next, the critic network was also taken to be an RBF neural network with six neurons in the hidden layer. The parameters of the controller and the adaptive critic networks are initialized randomly from a zero mean, isotropic Gaussian, with unit variance scaled by the fan-in of the output units. The target values of the critic network is then calculated using~(\ref{eq:probcriticNonL29}) for the specified range of $x_t$ and the critic training is carried out using scaled conjugate gradient method until the weights of the network have converged. The termination criterion of the training process is set to $\mid \bigtriangleup F(\chi) \mid<0.001$ and $\mid \bigtriangleup (\chi) \mid<0.001$ (both must be satisfied) within the limit of $10000$ iterations. Here $\mid \bigtriangleup F(\chi) \mid$  and $\mid \bigtriangleup (\chi) \mid$ is the absolute difference of the error function (defined as the sum of the squares of the errors between target and desired output values) and absolute difference of connection weights of the network between two successive iterations respectively. During training of the critic network, the weights of the action network are kept fixed. The output from the converged critic is then used in~(\ref{eq:probcriticNonL33}) solving for optimal control values and the action network is then trained for the same number of iterations, termination criterion and training method as that of the critic network. During training of the action network, the weights of the critic network are kept fixed. After the action network converged, the critic network is trained again (by adapting weights of the previously converged critic) using the outputs of the converged action network. The training of the critic and the action networks are alternated for $3$ cycles after which all adaptation is halted and the controller network's ability to return the plant state to the origin is tested.

The control quality of the controller designed in the above manner is then compared with the conventional DHP adaptive critic technique~\cite{Werbos92}. The same forward neural network model as that used in the above probabilistic design method is used in the conventional DHP critic to represent the forward dynamics of the system. However, only the deterministic forward model represented by the sum of two RBF networks, $\hat{x}_t = h(x_{t-1}) + g(x_{t-1}) u_t$ is required for implementation of the conventional DHP adaptive critic method. Moreover, For fair comparison, same noise sequence, initial conditions, and same structure of critic and control neural networks were used during the implementation of each control method.

\begin{figure}[htpb]
  \centerline{
    \begin{tabular}{c}
      {\includegraphics[width=3.0in]{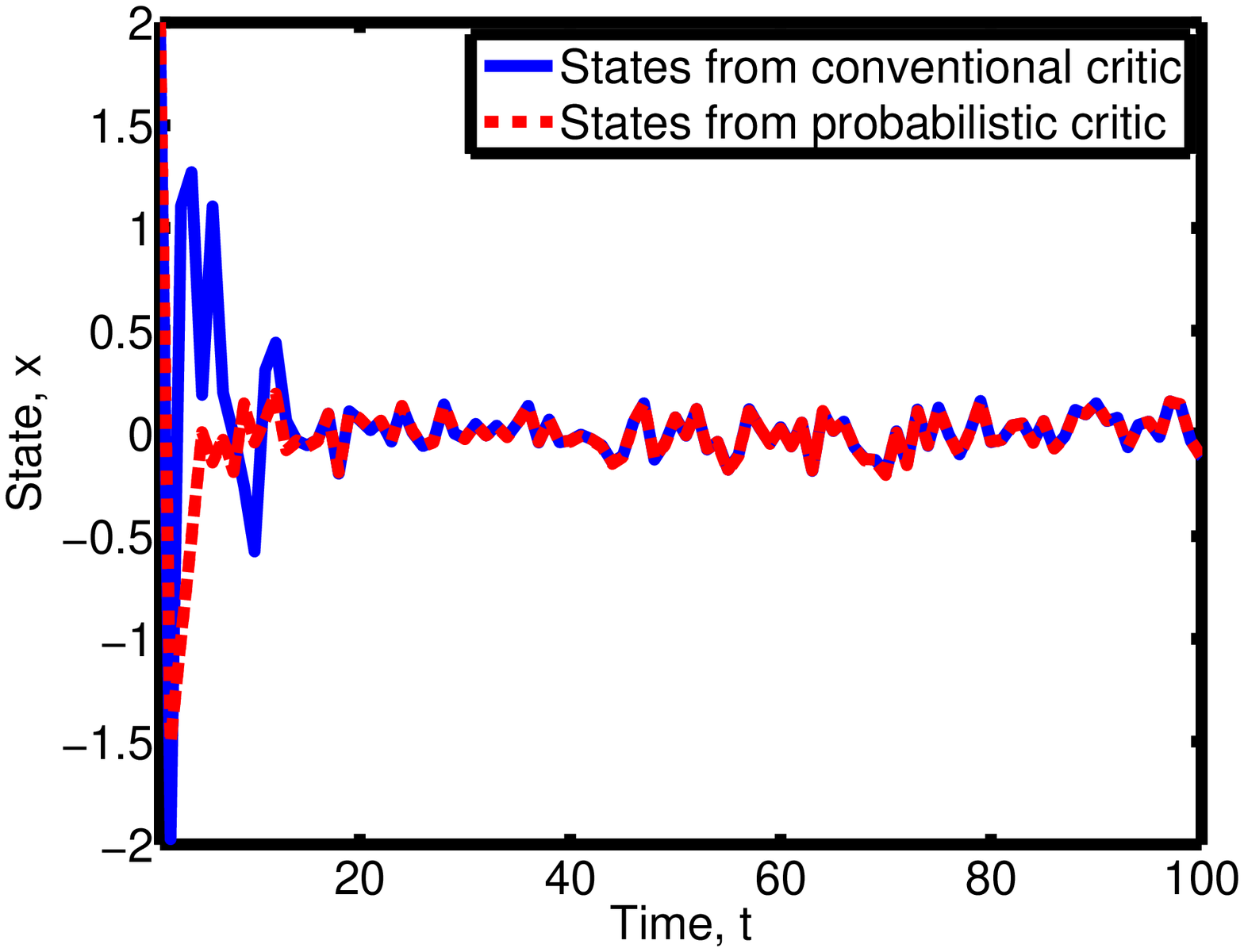}}
      \\
      (a)
      \\
      {\includegraphics[width=3.0in]{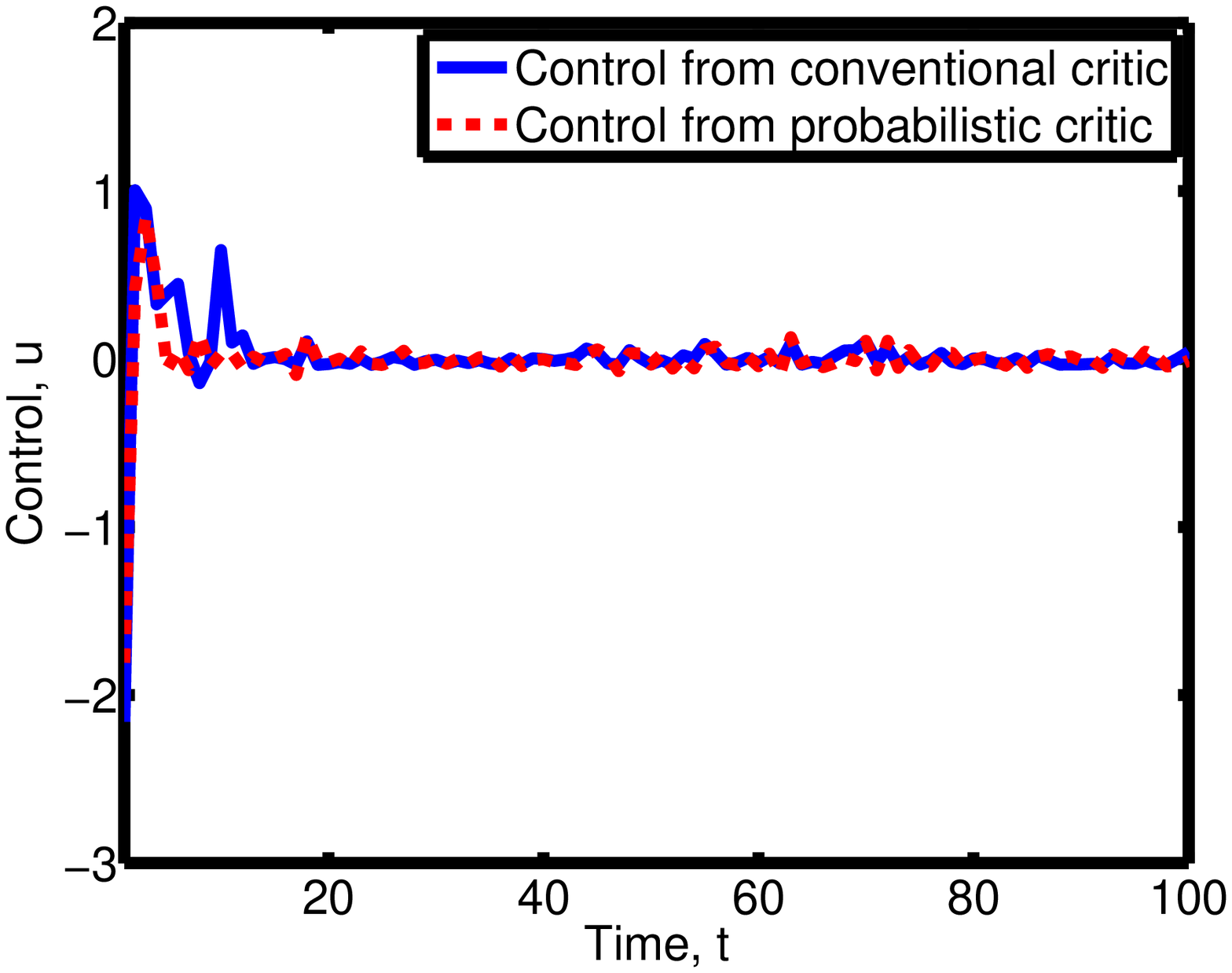}}
      \\
      (b)
    \end{tabular}
    }
  \caption[Nonlinear stochastic system.] {Nonlinear stochastic system: (a) Probabilistic and conventional critic estimated values for the state (b) Probabilistic and conventional critic values for control.}
  \label{fig:ProbCriticFig1}
\end{figure}

The ability of the obtained controllers from the conventional and probabilistic DHP adaptive critic methods in returning the system state from its initial value which is taken to be $x_0 = 2$ to the origin is then tested. Note that this initial condition was chosen arbitrarily. During testing, the optimal control signal, $u_t$ is generated from the controller network at each time instant $t$ based on previous state value $x_{t-1}$.  This optimal control signal is then forwarded to the plant of Equation~(\ref{eq:probcriticNonL34}) and the system output, $x_t$ is measured. Figure~\ref{fig:ProbCriticFig1} shows histories of the probabilistic and conventional DHP adaptive critic states and control efforts. Both of the probabilistic and conventional critic methods managed to regularize the state of the system around zero as required.  The probabilistic DHP critic method, however, ensures minimal overshoot compared to the conventional DHP critic method. This is expected since in the conventional critic method both the critic and controller network parameters have been tuned based on the assumption of the existence of an accurate deterministic forward model. In the probabilistic critic method on the other hand, optimal control law is derived such that the distance between the joint pdf of the closed loop system and an ideal joint pdf is minimized. This allows considering system's uncertainty and improves the performance of the derived optimal control law.

\section{Conclusion}\label{Conc}
In this paper, the solution to the probabilistic DHP adaptive critic method and randomized control input design have been addressed for a class of dynamic stochastic nonlinear systems. Using RBF neural networks to approximate the conditional expectations of the pdfs and unknown nonlinear models, a randomized control strategy has been developed which minimizes the discrepancy between the joint pdf of the closed loop system and a predetermined ideal joint pdf. It has been shown that the controller design is a nonlinear optimization problem and hence need to be solved numerically. A simulated example is used to illustrate the use of the randomized control algorithm as derived from the probabilistic critic and compared against its counter part of conventional critic design methods. The probabilistic design of critic networks showed minimum overshoot compared to the conventional critic design method.


\end{document}